\documentclass[STIX1COL]{WileyNJD-v2}
\usepackage{moreverb}
\usepackage{subfig}
\usepackage{bm}
\usepackage{amsthm}
\usepackage{amsmath}
\usepackage[T1]{fontenc} 
\usepackage{nicefrac} 
\usepackage{tikz,caption}
\usepackage{pgfplots}
\pgfplotsset{compat=newest}
\pgfplotsset{plot coordinates/math parser=false}
\usepackage{url}
\usepackage{color}
\usepackage{stix}

\usepackage{cleveref}
\newlength\figureheight
\newlength\figurewidth

\articletype{Article Type}%

\received{<day> <Month>, <year>}
\revised{<day> <Month>, <year>}
\accepted{<day> <Month>, <year>}

\newcommand\N{\ensuremath{\mathbb{N}}}
\newcommand\R{\ensuremath{\mathbb{R}}}
\newcommand\K{\ensuremath{\mathbb{K}}}

\newcommand\X{\ensuremath{{A}}}
\newcommand\x{\ensuremath{{a}}}

\newcommand{\abs}    [1] {\ensuremath{\left\lvert #1 \right\rvert}}
\newcommand{\norm}    [1] {\ensuremath{\left\lVert #1 \right\rVert}}

\newcommand*{\horzbar}{\rule[.5ex]{2.5ex}{0.5pt}}

\newcommand{\change}[1]{\textcolor{black}{#1}}

\begin{document}

\title{A literature survey of matrix methods for data science\protect\thanks{}}

\author[1]{Martin Stoll*}




\address[1]{\orgdiv{Chair of Scientific Computing, Department of Mathematics}, \orgname{TU Chemnitz}, \orgaddress{\state{Chemnitz}, \country{Germany}}}

\corres{Martin Stoll \email{martin.stoll@mathematik.tu-chemnitz.de}}

\presentaddress{Chair of Scientific Computing, Department of Mathematics, TU Chemnitz, Reichenhainer Str. 41, 09126 Chemnitz, Germany}

\abstract[Abstract]
{
	\change{Efficient numerical linear algebra is a core ingredient in many applications across almost all scientific and industrial disciplines. With this survey we want to illustrate that numerical linear algebra has played and is playing a crucial role in enabling and improving data science computations with many new developments being fueled by the availability of data and computing resources. We highlight the role of various different factorizations and the power of changing the representation of the data as well as discussing topics such as randomized algorithms, functions of matrices, and high-dimensional problems. We briefly touch upon the role of techniques from numerical linear algebra used within deep learning.}
}


\jnlcitation{\cname{%
		\author{}, 
	} (\cyear{2020}), 
	\ctitle{}, \cjournal{}, \cvol{}.}

\maketitle

\section{Introduction}\label{sec1} 
The study of extracting information from data has become crucial in many field ranging from business, engineering, fundamental research, or culture. \change{We here assume that data science intends to analyze and
	understand actual phenomena with \textit{data} according to 
	\cite{hayashi1998data}. To achieve this, we follow \cite{donoho201750} in that data science draws on elements of machine learning, data mining and many other mathematical fields such as optimization or statistics .} Also, we want to point out that in order to obtain information from data it is not necessarily implied that the amount of data is \textit{big} but often it is.

The multitude of applications where such data occur and need to be studied goes beyond the scope of this paper. Naturally, matrices arise as part of spatial data analysis \cite{han2001spatial,miller2009geographic,grilli2017review} or time-series data analysis \cite{aghabozorgi2015time,jeong2011weighted,liao2005clustering,rousseeuw2005robust,brown2004smoothing} but can also be found in many more disciplines \cite{van2016data}. 

We here view the data as being represented in a matrix 
\begin{equation}
\X=
\begin{bmatrix}
\horzbar&\x_1^T&\horzbar\\
&\vdots\\
\horzbar&\x_n^T&\horzbar
\end{bmatrix}
\in\R^{n,m}
\end{equation}
with the dimensions $m,n$ related to the underlying data. Here $m$ is the dimension of the feature space with feature vectors $\x_i\in \R^m$ viewed here as the rows of $\X$. The dimension $n$ is the number of data points and is thus possibly very large.  We often assume that $n>m$. Alternatively, the data can also arise in the form of a tensor of order $d$
\begin{equation}
\X\in \R^{n_1\times n_2\times \ldots\times n_d}
\end{equation}
with dimensions $n_i\in\N\quad \forall i$.

The tasks of extracting meaningful information from the collected data varies between application areas. \change{ In the process of extracting information from data we often encounter tasks such as \emph{unsupervised learning}, \emph{semi-supervised learning}, and \emph{supervised learning} (cf.  \cite{James_2013,Hastie_2009}). The difference between these can roughly be summarized by the availability and usage of training data with none for unsupervised, only for a subset for semi-supervised, and for all of the data in supervised learning.}

Before starting the detailed discussion, we would like to point out a particular example that is a core problem in data science, statistics, numerical linear algebra and computer science alike. This is the least squares problem \cite{GoluVanl96}, where in brief one wants to fit a linear model to labeled data. Let us assume that we are given $n$ data points $\x_i\in\R^m$ and typically outputs $y_i$, e.g. a value of either $-1$ or $1$ for a classification problem. We are then interested in finding a weight vector $w\in\R^m$ such that the $2$-norm of the residual 
$
r_i=\x_i^Tw-y_{i}
$
is minimized.It is well known that this would lead to a linear regression problem
\begin{equation}
\label{ls1}
\min_{w}\norm{\X w-y}_2^2+\lambda \norm{w}_{2}^2
\end{equation}
with 
$\X$ as above.
Here,  we also introduce $\lambda>0$ as a regularization or ridge parameter for the regularization term $\norm{w}_{2}^2.$ The regularization term could also be measured in several different norms such as the $l_1$-norm \cite{tibshirani1996regression} or the total variation norm \cite{vogel1996iterative}. The solution to this problem is then given by 
\begin{equation}
\label{ls2}
w_{*}=\left(\X^{\top}\X+\lambda I\right)^{-1}\X^T y,
\end{equation}
a prototypical linear system of equations with a symmetric and positive definite matrix as long as $\lambda>0.$ \change{This example illustrates that the different disciplines are intertwined. Problem \eqref{ls1}  arises in data mining relying on techniques from numerical linear algebra to make the evaluation robust and efficient. On the other hand, the development of sophisticated numerical methods is driven by studying problems with \textit{ real data}.} The goal of this survey is to show how information extraction from data, data modeling, and pattern finding relies on efficient techniques from numerical linear algebra that not only enable computations but also reveal hidden information. 

The paper is structured as follows. We first illustrate the use of classical factorizations such as the QR and singular value decomposition (SVD) and then introduce interpretable factorizations such as the non-negative matrix factorization (NMF) or the CUR decomposition. We additionally discuss literature devoted to kernel methods with special attention given to the graph Laplacian. Randomization as well as functions of matrices are discussed next. We close with a brief discussion of recent results for high-dimensional problems and deep learning applications.

As a word of caution, we want to remark that the field of numerical linear algebra is vast and the analysis of data via techniques from numerical linear algebra is not new. The goal of this survey is to point to recent trends and we apologize to the authors whose results we missed while writing this. In particular we want to refer to the beautiful books by Eld{\'e}n \cite{elden2019matrix} and Strang \cite{strang2019linear} that provide general introductions to linear algebra for data science applications.
\section{Data matrices and factorizations}
\label{sec::fac}
The decompositional approach to matrix computations has been named one of the top 10 algorithms of the 20th century \cite{dongarra2000guest}.  Matrix factorizations are an ubiquitous tool in data science and have received much attention over the last years. A great example is the use of matrix factorization techniques for recommender systems such as the Netflix challenge\footnote{\url{https://en.wikipedia.org/wiki/Netflix_Prize}} \cite{koren2009matrix}. We here review some important matrix factorizations, their applications as well as tailored factorizations used for data science.
We split the discussion into classical factorizations, which have been the workhorse of many applications such as engineering or fluid mechanics, and factorizations that are designed to more closely resemble the nature of the data.
\subsection{Classical factorizations}
\subsubsection*{The singular value decomposition}
Given a data matrix $\X\in\R^{n,m},$ assuming $n\geq m$, a singular value decomposition (SVD) \cite{GoluVanl96} is given as
\[
\X=USV^{\top}
\]
with $U\in\R^{n,n}$ and $V\in\R^{m,m}$ being orthogonal matrices. The matrix $S\in\R^{n,m}$ is of the following form
\[
S=
\begin{bmatrix}
\sigma_1&0&\ldots&0\\
0&\sigma_2&&\\
&&\ddots&\\
\vdots&&&\sigma_m\\
0&\ldots&\ldots&0\\
\vdots&\vdots&\vdots&0\\
0&\ldots&\ldots&0\\
\end{bmatrix}
\]
with singular values $\sigma_1\geq\sigma_2\ldots\geq\sigma_m\geq 0.$ An equivalent and often very useful representation is the outer product form of the SVD as
\[
\X=\sum_{i=1}^{m}\sigma_iu_{i}v_{i}^{T}
\]
where $u_i$ and $v_i$ are the columns of $U$ and $V$, respectively.
\change{
This allows for a natural interpretation of $U$ as providing a basis for the column space of $\X$ and $V$ for its row-space. A crucial task is to find a good rank-$k$ approximation to the matrix $\X$.  The Eckart--Young--Mirsky theorem \cite{GolVLoan13} states that the best rank-$k$ approximation in any unitarily invariant matrix norm is given as
\[
\X\approx \X_k:=\sum_{i=1}^{k}\sigma_iu_{i}v_{i}^{T},
\]
hereby ignoring the smaller singular values in the summation. This is known as the \emph{truncated singular value decomposition} and is written in matrix form as
\begin{equation}
\label{tsvd}
\X_k=U_kS_{k}V_{k}^{\top}.
\end{equation}
}
The SVD is one of the most crucial tools in the complexity reduction of large-scale problems. The singular value decomposition is naturally well-suited to solve the least squares problem \eqref{ls1} of the form 
\begin{equation}
w_{*}=\left(\X^{\top}\X+\lambda I\right)^{-1}\X^{\top} y=
\left(VS^{\top}U^{T}USV^{\top}+\lambda I\right)^{-1}VS^{\top}U^{\top}y=
V\left(S^{\top}S+\lambda I\right)^{-1}S^{\top}U^{\top}y.
\end{equation}
This is expensive due to the cost of computing the full SVD but the truncated SVD can be exploited for solving least squares problems as discussed in \cite{hansen1987truncatedsvd}.  There have been many algorithmic updates in the computation of the (truncated) SVD that deal with the numerical difficulties of large-scale problems, rounding errors, locking of wanted singular vectors and purging of unwanted information \cite{baglama2005augmented,hochstenbach2004harmonic,stoll2012krylov,andez2008robust}.  At the heart often lies the Lanczos bidiagonalization and for large scale problems incorporating implicit restarts is mandatory (cf. \cite{bjorck2014stability, hernandez2007restarted, baglama2005augmented, stoll2012krylov}). \change{The algoritmic foundations for computing the truncated SVD, namely the Lanczos bidiagonalization \cite{golub1965calculating} was shown to be equivalent \cite{elden2004partial,bjorck2014stability} to a well-known method in statistics, namely the \textit{NIPALS} (Nonlinear Iterative Partial Least Squares) method.} To the best of our knowledge, the algorithmic improvements developed for the Lanczos-bidiagonalization have not been  exploited within NIPALS implementations (cf. \cite{bjorck2014stability}). Nevertheless, NIPALS has become a crucial tool in the analysis of problems from economics applications \cite{hulland1999use,hair2016primer,f2014partial}. 

The SVD is also a key ingredient in model order reduction\footnote{A prototypical example is the drastic reduction of the dimensionality of the system matrices defining a dynamical system.} \cite{benner2015survey,chaturantabut2010nonlinear} classically used for reducing the dimensionality of physics-related models based on differential equations. Recently, model order reduction has found more and more applications in machine learning \cite{kutz2017deep,trehan2017error}. One of the most important applications that directly mirrors the use of the SVD in computational science and engineering is the creation of reduced representations. Here, applications include text mining \cite{albright2004taming}, face recognition \cite{zhang2010discriminative}, medicine \cite{zear2018proposed} and many more.

The SVD also comes in many disguises among the different disciplines of statistics, engineering, and applied mathematics. Given a data matrix $\X$ applying principal component analysis (PCA), which is equivalent to performing the SVD, has been a key tool for understanding the structure of the data. In PCA, the column means within $\X$ is zero and then the right singular vectors $v_i$ are called the \textit{principal component directions} of $\X$ and the left singular vectors are the \textit{principal components of $\X$}. More information and applications are given in \cite{jolliffe2011principal,jolliffe2016principal,metsalu2015clustvis,yi2017joint}.

In order to compute the singular value decomposition for data matrices that originate from massive datasets one often has to resort to techniques from high performance computing \cite{anderson1990lapack,sukkari2016high,nakatsukasa2013stable}. Additionally, it is possible to rely on randomized algorithms that we discuss in Section \ref{sec::random}.

The SVD is also a crucial ingredient in many algorithms for high-dimensional data analysis, see Section \ref{sec::tensors} on tensor factorizations.

\subsubsection*{The QR factorization}
Given a set of vectors collected in the matrix $A$ and considering the span of the columns of $A$ it is clear that the vectors themselves can potentially \change{be a terribly conditioned basis for further numerical computations.} Hence, one wants to find a well-conditioned basis, ideally a basis with orthogonal vectors. This task is achieved by computing the QR factorization, i.e., 
\[
\X=QR
\]
with $Q\in\R^{n,n}$ an orthogonal matrix and $R\in\R^{n,m}$ an upper triangular matrix of the form
\[
R=
\begin{bmatrix}
\hat{R}\\
0
\end{bmatrix},
\]
where the matrix $\hat{R}\in\R^{m,m}$ is invertible if only if the matrix $\X$ has full column rank $m$. From the representation it is clear that one can also work with 
$
\X=\hat{Q}\hat{R}
$
where $\hat{Q}$ only contains the first $m$ columns of $Q$. The QR factorization is another ubiquitous factorization in applied mathematics. Its computation is typically rather expensive and for the details we refer to \cite{GoluVanl96}. In particular, the reduction of $\X$ to triangular form via so-called Householder reflectors is a de-facto standard. The solution of the least squares problem \eqref{ls1} without regularization \cite{GoluVanl96} via the QR factorization is well-known 
$$
w_{*}=\mathrm{argmin}_w\norm{QRw- y}_2=\mathrm{argmin}_w\norm{Rw-Q^\top y}_2.
$$
%
Truncated QR decompositions in disguise are also at the heart of the Lanczos \cite{lanczos1950iteration} and Arnoldi \cite{arnoldi1951principle} methods, which are the key algorithms of Krylov subspace methods. \change{Many algorithms in numerical linear algebra rely on Krylov subspaces, i.e.
	$$
	\mathcal{K}_\ell(M,r):=\mathrm{span}\left\lbrace r,Mr,M^2r,M^3r,\ldots,M^{\ell-1}r\right\rbrace
	$$
	being the space of dimension $\ell$, where $M$ is the system matrix of the underlying problem, e.g. $M=\left(\X^{\top}\X+\lambda I\right)$ for the ridge regression problem \eqref{ls2}, and $r$ is a vector associated with a right-hand side of \eqref{ls2}. In order to obtain robust methods we rely on a well-conditioned basis of $\mathcal{K}_\ell(M,r)$. The Lanczos method computes an orthonormal basis of increasing dimensionality with great efficiency only requiring one matrix vector product per iteration. For more details, we refer to \cite{GolVLoan13,saad2003iterative}.
}
Truncated QR decompositions are based on computing the columns of a matrix $Q$ spanning the basis of the range of the Krylov matrix. 

Additionally, the pivoted QR factorization \cite{stewart1998matrix,stewart1999four} computes the factorization
$$
\X P=QR
$$
with $P$ a permutation matrix. This QR factorization can be computed maintaining the sparsity of the matrix $\X$ (cf. \cite{stewart1999four}) and it is also closely related to interpretable factorizations, which we discuss next.

\subsection{Interpretable Factorizations}
\label{sec::ID}
Despite the beautiful mathematical properties of both the QR and the SVD they sometimes do not provide a representation that allows an easy interpretation of the results for practitioners. \change{For example, in many applications the data are nonnegative while the singular vectors can contain negative values.} Mathematically this means that a given data matrix $\X\in\R^{n,m}$ is well approximated by the truncated SVD  $\X\approx U_kS_{k}V_{k}^{\top}$ but the singular vectors contained in $U_k(:,l),$ are in general not sparse or non-negative even though this holds for the data. \change{We adopt \textsc{Matlab} notation addressing columns and rows of matrices by using $\X(:,k)$ for the $k$-th column of $\X$ and $\X(k,:)$ for the $k$-th row. Analogously, this can be defined for several rows or columns.}

\begin{figure}[htb!]
	\begin{center}
		\setlength\figureheight{0.4\linewidth} 
		\setlength\figurewidth{0.6\linewidth}
		\includegraphics[width=8cm]{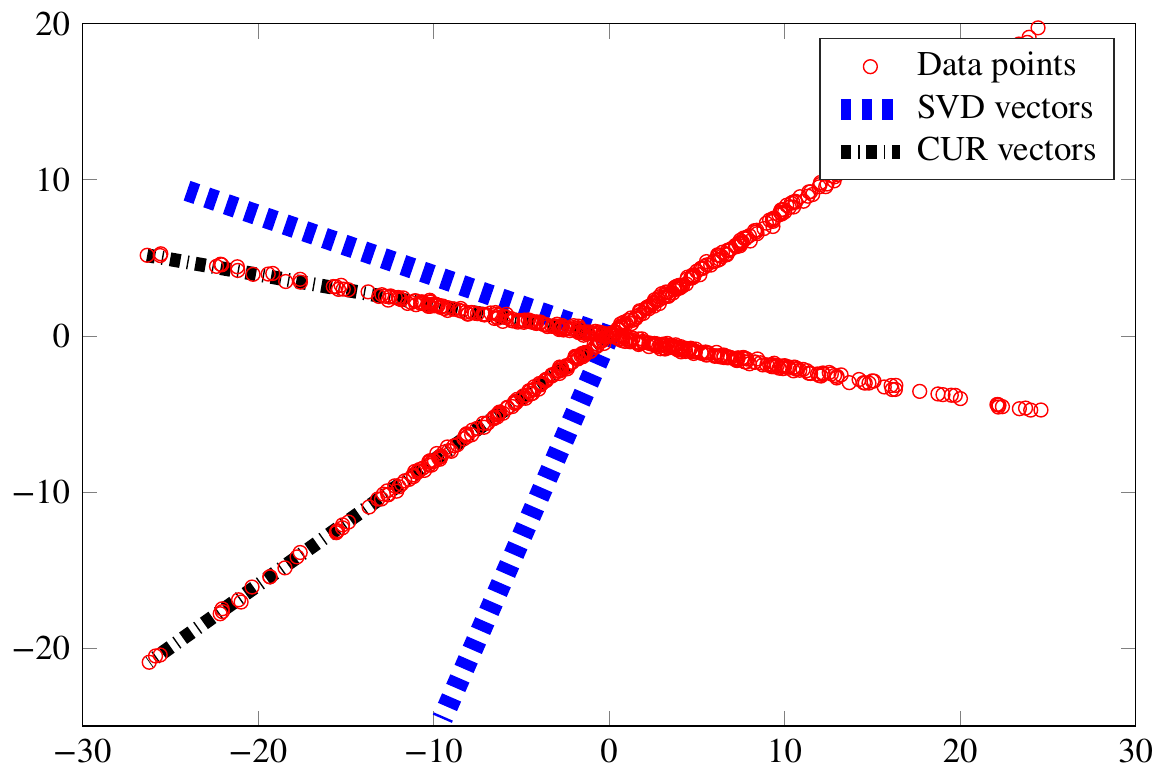}
		\caption{We here illustrate a set of datapoints (red) that are stored in the matrix $\X$. We then show the two dominating rows (black) of $\X$ obtained from a CUR decomposition and the two dominating singular vectors (blue).}
		\label{fig:CUR}
	\end{center}
\end{figure}

To preserve the properties inherent in the application while maintaining a good approximation of the original data with reduced complexity, many interpretable factorizations have been developed over recent years (cf. \cite{mahoney2009cur,woodruff2014sketching,drineas2012fast,boutsidis2017optimal,lee2001algorithms,donoho2004does,gillis2014and} for some of them). 
Here, \cref{fig:CUR} illustrates how the CUR approximation, which we introduce later in this section, naturally represents the original data.

We now briefly introduce these methods and comment on applications and new developments. We here follow the notation of \cite{voronin2017efficient} and start with the CX decomposition defined by
\[
\X\approx CX
\]
where $C\in\R^{n,k}$ and $X\in\R^{k,m}$. Here the matrix $C$ is a matrix consisting of $k$ columns of the data matrix $\X$, i.e., 
\[
C=\X(:,J)
\] for an index set $J$ of order $k$. The important feature of the matrix $C$ is that its columns are interpretable as they are taken from the original data. Since the data matrix is often sparse this is inherited by $C$ and  as a result storing $C$ requires less memory than storing the singular vector matrix $U_k$. The computation of $C$ and $X$ is done via the minimization of 
\[
\min\norm{\X -CX},
\]
where $\|\cdot\|$ may be the 2-norm or the Frobenius norm. This problem is known as the column subset selection problem \cite{boutsidis2009improved} and its NP-completeness is discussed in 
\cite{shitov2017column}. Typically, more structure is required than just a general matrix $X$ and one quickly moves to the interpolative decomposition (ID) \cite{voronin2017efficient}
\begin{equation}
\X\approx CV^{\top}
\label{IDr}
\end{equation}
where $C\in\R^{n,k}$ is as before but the matrix $V\in\R^{m,k}$ is constructed such that it contains a $k\times k$ identity matrix and $\max_{i,j}\abs{v_{ij}}\leq 1$. A procedure to obtain a low-rank ID via a pivoted QR is given in \cite{voronin2017efficient} returning a $V^T=\left[I_k\  T_l\right]P^T$ with $P$ a permutation matrix and $T_l$ a solution related to the upper triangular factors of the pivoted QR.  This approach works analogously when a one-sided representation in terms of matrix rows is desired. 
\change{
One can also obtain a two-sided interpolative decomposition 
\begin{equation}
\label{2ID}
\X\approx W\X(I,J)V^T.
\end{equation}
We start its computation by using a one-sided interpolative decomposition $\X\approx CV^{\top},$ which already provides us with the set of crucial column indices $J$. In order to compute the row indices $I$ and the matrix $W$, we now compute a one-sided decomposition for the matrix $C^T$. Following \eqref{IDr} we then obtain
$$
C^T\approx \tilde{C}\tilde{V}^T,
$$
where $\tilde{C}$ contains by design columns of $C^T$, i.e., rows of $C$ and thus elements of the rows of $\X$ along with the row index set $I$. The matrix $W$ is then given as $W=\tilde{V}$.}
Note that $W$ and $V$ do not contain columns and rows of the original matrix $\X$, respectively, and hence properties such as sparsity and non-negativity found in the data are typically not carried over.

\change{
As a result, we can consider a factorization that avoids this pitfall and bears a lot of similarity to the two-sided ID.  This is achieved by the CUR decomposition
\[
\X\approx CUR
\]
where $C\in\R^{n,k},$ $U\in\R^{k,k},$ and $R\in\R^{k,m}.$ Here, $C$ contains columns of the original matrix and $R$ represents a subset of its columns. Both matrices inherit properties such as non-negativity, sparsity, and finally interpretability. }The CUR decomposition is also known as the skeleton decomposition \cite{kishore2017literature,gantmacher1964theory,tyrtyshnikov1996mosaic,osinsky2018pseudo,goreinov1997pseudo} and is closely related
to a rank-revealing QR factorization \cite{voronin2017efficient,gu1996efficient}. 
\change{
 The point of departure for computing the CUR decomposition is typically a low-rank factorization of the matrix $\X$. Both the truncated SVD \eqref{tsvd} and the two-sided ID \eqref{2ID} are the typical initial factorizations.
 Given a truncated SVD, the crucial algorithmic step is to select the index sets $I$ and $J$. For this one typically computes \textit{leverage scores} as sums over the rows of the matrices $U_k$ and $V_k$ coming from the truncated SVD
$\X\approx U_kS_kV_k^{\top},$ e.g. $\ell_j=\sum_{i=1}^{k}U_{j,i}$ for $j=1,\ldots,n$ for the column selection and analogously for $V_k$. In \cite{mahoney2009cur} the authors provide a statistical interpretation of the leverage scores as a probability distribution. More recently, Embree and Sorensen \cite{sorensen2016deim} have introduced a procedure based on the discrete empirical interpolation method (DEIM) \cite{chaturantabut2010nonlinear} where the selection of the column and row indices is based on a greedy projection technique resembling a pivoting strategy within the LU factorization. In \cite{deshpande2010efficient} the authors compute the column and row subset using conditional expectations with a more efficient numerical realization being recently introduced in \cite{cortinovis2019low}.
It remains to compute the intersection matrix $U$ as $U=\X(I,J)^{-1}$ with the other possibility being $U=C^{\dagger}AR^{\dagger},$ where $\dagger$ indicates the Moore-Penrose inverse. Recently, a perturbation analysis of the CUR decomposition was presented in \cite{hamm2019perturbations}.}

Another important interpretable factorization of the matrix $\X$ is the so-called non-negative matrix factorization (NMF) \cite{Gillis2014,elden2019matrix,ding2008convex,ding2008convex,berry2007algorithms}
\begin{equation}
\X^{T}\approx WH,
\end{equation}
which is a low-rank approximation using $W\in\R^{m,k}$ and $H\in\R^{k,n}$ with component-wise non-negativity written as $W\geq0,~H\geq 0.$ The interpretation of the columns of $W(:,j)\in\R^{m}$ is that they form a basis of order $k$ that best approximates the data points $\x_j,$ i.e., the columns of $\X^{T}$ via
$$
\x_l=WH(:,l)=\sum_{j=1}^{k}W_{:,j}H_{j,l}.
$$ 
The coefficients of how the data are expanded in the basis defined by $W$ are stored in $H$. Such linear dimension reduction frameworks are found in various data tasks within image processing or text mining. Here again the non-negativity of the elements in $W$ means that the resulting matrix is more sparse than an SVD-based approach and provides an interpretable feature representation \cite{liu2011constrained,li2016graph}.The computation of a NMF is an NP-hard ill--posed problem \cite{vavasis2009complexity} and it is typically based on solving the minimization problem 
\[
\min_{W,H\geq 0}\norm{\X^{T}-WH}_F.
\]
Alternating minimization procedures, which consist of an alternating update of the factors $W$ and $H$,  are typically employed and we refer to \cite{gillis2014and,hajinezhad2016nonnegative,luo2015nonnegative} for more details. \change{It is also possible to include further constraints such as sparsity as was done in \cite{rapin2014nmf,rapin2013sparse}.} In \cite{ding2005equivalence} the authors analyze the relationship between the NMF factorization of a \emph{kernel matrix} and spectral clustering discussed later. For further improving the performance the authors in \cite{cai2010graph} include a kernel matrix as a regularization term for the objective function. This shows that kernel matrices, which are matrices changing the representation of the data, are often very useful and we discuss these next.
%
\section{Changing the data representation} 
So far we focused on methods that directly utilize the data $\X$ as a matrix. Often it is necessary to transform the data to a different representation. The goal is that for the transformed data the learning task is easier and we now describe several approaches designed for that purpose.
\subsection{The graph Laplacian operator}
\label{sec::GL}
The data encoded in 
$\X\in\R^{n,m}$
either have a natural representation as a graph with the nodes $v_j$ representing the associated feature vectors $\x_j\ \forall j$ or they can be modeled that way. The result is a graph $G=(V,E)$ consisting of the nodes $v_j\in V$ and edges $e\in E,$ where an edge $e$ consists of a pair of nodes. We here consider undirected graphs where an edge is typically equipped with an edge weight representing the \textit{strength} of the connection between the corresponding nodes. Practically, the most relevant weight function is the Gaussian weight function 
\begin{equation}\label{eq::Gaussian}
w(v_i,v_j)=w_{ij}=\exp{\left(\nicefrac{-\norm{\x_i-\x_j}_2^2}{\sigma^2}\right)}.
\end{equation}
\change{
Many applications naturally have a graph structure but one can also convert data to graph form and we refer to \cite[Section 2.2]{von2007tutorial} where different techniques are presented.  Given a graph the weights are collected into a matrix $W\in\R^{n,n}$ where the diagonal is set to zero. If two nodes are not connected in the graph the associated matrix entry is set to zero. A particularly interesting and challenging example is a fully connected graph where all data points are compared pairwise also resulting in a dense matrix $W$.} As a second ingredient we compute the diagonal degree matrix $D$ where $d_{ii}=\sum_{j=1}^{n}w(v_i,v_j).$  \change{We then obtain the \emph{graph Laplacian} 
$L=D-W,$ which is often used in a normalized form either as the symmetric normalized Laplacian $L_{sym}=I-D^{-1/2}WD^{-1/2}$ or as the random walk Laplacian $L_{rw}=D^{-1}L$.
The properties of the graph Laplacian are discussed in \cite{von2007tutorial,chung1997spectral}.
In more detail, we can see that $L$ is symmetric and its positive semi-definiteness follows from
\begin{equation}
u^{\top}Lu=\frac{1}{2}\sum_{i,j=1}^{n}w_{ij}\left(u_i-u_j\right)^2.
\end{equation}
This relation only changes slightly when $L_{sym}$ is used.} 
In the context of data science many of the properties of the graph Laplacian can be utilized for tasks such as clustering. In particular, the eigeninformation of $L$ and $L_{sym}$ are crucial. For example,  the number of zero eigenvalues gives information about the number of connected components of the underlying graph. The eigenvector corresponding to the first non-zero eigenvalue is known as the {Fiedler vector}. As the eigenvector corresponding to the zero eigenvalue is the constant vector $c\begin{bmatrix}1,\ldots,1\end{bmatrix}^T$  with $c\in\R$ and since the Fiedler vector is orthogonal to it, we must have sign changes in the Fiedler vector. This makes the Fiedler a first candidate to perform clustering simply by the sign of its entries\footnote{In \url{https://people.eecs.berkeley.edu/~demmel/cs267/lecture20/lecture20.html} a connection to vibrating strings and standing waves is made that beautifully illustrates this property.} \cite{urschel2018nodal,ding2001min,bertrand2013seeing,chan1997optimality}.

\begin{figure}[htb!]
	\centering
	\subfloat{\includegraphics[width=.45\linewidth]{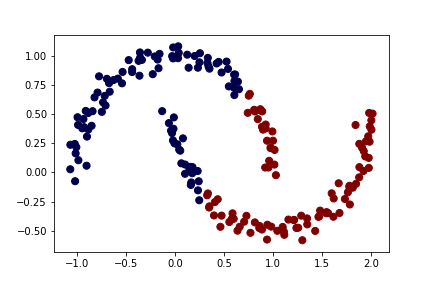}}
	\hfill
	\subfloat{\includegraphics[width=.45\linewidth]{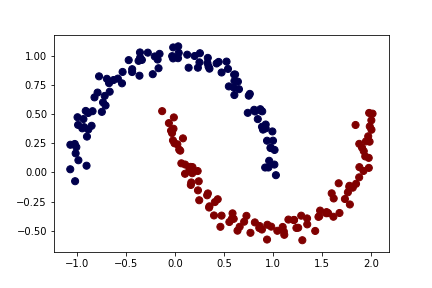}}
	\caption{The typical \textit{two-moons} data set with desired unsupervised classification into two-classes. A standard \textit{scikit-learn}\protect\footnotemark $k$-means clustering applied on the left  vs. \textit{scikit-learn} spectral clustering (right). The curvature of the data is only correctly identified with the spectral clustering approach.}
	\label{fig:clustering}
\end{figure}

In fact, one of the classical tasks that is performed using the eigeninformation of the graph Laplacian is spectral clustering \cite{von2007tutorial}, where one computes the first\footnote{The eigenvalues of the Laplacian are given as $0 = \lambda_1 \le \lambda_2 \le \dots \le \lambda_n$.} $k$ eigenvectors $\phi_1, \ldots,\phi_k$. As the graph Laplacian translates the original data $\X$ into an alternative space encoded into new matrices $W$ and $D$, its eigeninformation will lead to a different clustering behavior than traditional methods such as {k-means} \cite{steinley2006k}. It can be seen from \cref{fig:clustering} that standard k-means clustering based on $\X$ \cite{hartigan1979algorithm,kanungo2002efficient} shows poorer performance when compared against spectral clustering \cite{von2007tutorial}. In more detail, the most common spectral clustering methods proceed by using k-means on the rows of the eigenvector matrix 
\begin{equation}
\label{phi}
\Phi_k=\left[\phi_1,\ldots,\phi_k \right].
\end{equation}
In \cite{von2007tutorial} the author illustrates that performing spectral clustering solves relaxed versions of known graph cut problems, which aim at partitioning the vertices of a graph into disjoint subsets. 
\footnotetext{\url{https://scikit-learn.org}}
If instead of the graph Laplacian the two normalized versions are used, one obtains different results, for $L_{rw}$ see \cite{shi2000normalized} and for $L_{sym}$ \cite{ng2002spectral}. The hyperparameter\footnote{As hyperparameter we understand a parameter whose value is set before the learning process starts.} $\sigma$ in the denominator of the weight function can be replaced by a local scaling with an additional hyperparameter describing the locality \cite{zelnik2005self}. Often the parameter is chosen in a heuristic way according to the performance of the algorithm using the graph Laplacian  \cite{ng2002spectral}. When interpreted in terms of Gaussian processes the parameter $\sigma$ is related to the characteristic length-scale of the process \cite{rasmussen2003gaussian}.

As the foundation of the spectral clustering method is the computation of $k$ eigenvectors of the (normalized) graph Laplacian, it is important to be able to compute these eigenvectors efficiently. Be reminded that the matrices $L,$ $L_{rw},$ and $ L_{sym}=I-D^{-1/2}WD^{-1/2}$ are singular and the multiplicity of the zero-eigenvalue corresponds to the number of connected components in the graph. We are interested in computing the smallest eigenvalues of $L$. Iterative eigenvalue algorithms typically converge towards the largest eigenvalues and as a result we would need to invert the matrix $L$, which is not possible. Focusing on $L_{sym}$ it is obvious that the smallest eigenvalues of $L_{sym}$ can be computed from the largest eigenvalues of 
$D^{-1/2}WD^{-1/2}.$ Methods for efficiently computing the eigeninformation of such a matrix often rely on Krylov subspaces \cite{GolVLoan13,stewart2001matrix,lehoucq1998arpack,saad2011numerical}, where it is crucial to perform the matrix vector products with $D^{-1/2}WD^{-1/2}$ efficiently. If we assume that the graph consists of one connected component $D$ is diagonal and invertible. The main cost of multiplying with $D^{-1/2}WD^{-1/2}$ comes from computing matrix vector products with $W$. For sparse graphs the matrix $W$ will itself be sparse and matrix vector products will be inexpensive. The main computational challenge is then encountered for non-sparse or fully connected graphs. For large data-sets matrix vector products with $W$ are often infeasible and hence more sophisticated techniques are needed. Recently, methods based on the non-equispaced Fourier transform \cite{alfke2018nfft}, the (improved) fast Gauss transform \cite{yang2005efficient,morariu2009automatic}, or algebraic fast multipole methods \cite{yu2017geometry,march2015kernel} have shown great potential resulting in a complexity of  $\mathcal{O}(n\log n)$ for the matrix vector products with a fixed number of columns $m$. While these methods provide great speed-ups the dimensionality $m$ of the feature vectors still provides a significant challenge in computations and we return to this point in the next section. 

The basis  $\Phi$ given in \eqref{phi} is not only important for spectral clustering but also for a group of methods recently introduced for
semi-supervised learning, i.e., methods that, given a small set of labeled data, classify the remaining unlabeled points simultaneously. The methods are based on partial differential equation techniques from material science modeling \cite{cahn1958free,allen1979microscopic}, namely diffuse-interface methods, and have also been used in image inpainting \cite{bertozzi2006inpainting}. The graph Laplacian then replaces the classical Laplacian operator resulting in a differential equation based on the graph data via
\begin{equation}
\label{ACgraph}
u_t=\varepsilon  L_{sym}u-\frac{1}{\varepsilon}\psi'(u)+\omega (f-u)
\end{equation}
\change{
with $\psi(u)$ being a potential defined on the graph nodes enforcing two classes and $\omega$ incorporating penalization for deviation from the training data stored in $f$.} The variable $\varepsilon$ is a hyperparameter related to the thickness of the interface region. Due to the large number of vertices the dimensionality of \eqref{ACgraph} is vast. A projection using $\Phi_k$ reduces the PDE to a $k$-dimensional equation. Similar to model order reduction the nonlinearity still needs to be evaluated in the large-dimensional space \cite{chaturantabut2010nonlinear}. 

We briefly want to comment on the use of the graph Laplacian in image processing, where the pixels are often represented as the nodes in a graph, which would then lead to a fully connected graph. In this field the graph Laplacian is often used as a regularizer for denoising \cite{kheradmand2013general,liu2014progressive,romano2017little} or image restoration \cite{kheradmand2014general}. In \cite{milanfar2012tour} the construction of the Laplacian is performed patch-wise in both a local and a non-local fashion. A more in-depth discussion for applications in image processing can also be found in \cite{cheung2018graph}.

The graph Laplacian has also recently enjoyed wide applicability within deep learning, namely, as an essential ingredient within so-called \emph{Graph Convolutional Networks} \cite{henaff2015deep,bruna2013spectral,kipf2016semi} where the equation at layer $l$ for semi-supervised learning becomes
\[ X^{(l+1)} = \sigma_l \left( \sum_{k=1}^{N_K} K^{(k)} X^{(l)} \Theta^{(l,k)} \right), \] 
with $\sigma_l$ an activation function and weights $\Theta$ that need to be learned. \change{The crucial filter matrices $K^{(k)}$ are computed using the eigeninformation of the graph Laplacian associated with the input $X^{(0)}$. }
The matrices $K^{(k)}$ are composed from a $k$-dimensional filter space $\mathrm{span}\left\lbrace\phi_1,\ldots,\phi_k\right\rbrace,$ typically of the form  $K^{(k)}=U\phi_k(\Lambda)U^T,$ where $U$ and $\Lambda$ are the eigenvector and eigenvalue matrix of the graph Laplacian or a slight modification of it. \change{The name convolutional network stems from the fact that the transformation $U\phi_k(\Lambda)U^T$ can be interpreted as graph Fourier transform \cite{shuman2012emerging,defferrard2016convolutional}. Similarly, classical convolutional networks apply a filter/convolution to the data to detect more structure within the data \cite{lecun1995convolutional}.} In more detail, for a signal $x\in\R^n,$ with $n$ the number of nodes in the graph, $\hat{x}=U^Tx$ is the graph Fourier transform and its inverse is given by $x=U\hat{x}$. It is clear that polynomial filters $\phi_j$ are easy to apply either directly, by multiplying with the matrices, or via an (approximate) eigendecomposition of the Laplacian. Many filters and efficient methods for their computations have been suggested and we refer to \cite{kipf2016semi,henaff2015deep,alfke2019semi,levie2018cayleynets,shuman2018distributed,li2019scalable,zhou2018graph,wu2019comprehensive} for some of them.

There have also been generalizations of the graph Laplacian to other settings. We in particular want to mention the case of \emph{hypergraphs} \cite{zhou2007learning,zhou2005beyond}, where an edge is now a collection of possibly many nodes. Again, one can obtain a normalized Laplacian operator of the form $L=I-\Xi$ and perform spectral clustering based on the eigeninformation of this matrix. Hypergraphs are encountered in many applications such as biological networks \cite{klamt2009hypergraphs}, image processing \cite{yu2012adaptive}, social networks \cite{zhang2010hypergraph} or music recommendation \cite{bu2010music}. Hypergraphs have also been used in the context of semi-supervised learning \cite{leordeanu2011semi,hein2013total,bosch2018generalizing} and particular in convolutional neural networks based on hypergraphs \cite{alfke2019semi,yadati2018hypergcn}.

\begin{figure}[htb!]
	\begin{center}
\includegraphics[width=8cm]{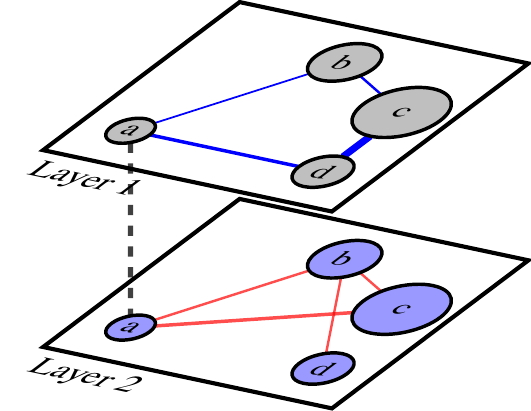}		
		\caption{A simple multilayer graph.\label{fig::mlgraph}}
	\end{center}
\end{figure}

Graph Laplacians have also been used to analyze multilayer networks \cite{kivela2014multilayer,boccaletti2014structure} where a set of nodes can be connected in various layers (cf. \cref{fig::mlgraph} for an illustration). The connections between the nodes and the various layers can be represented as a tensor but also using a (supra)-Laplacian \cite{sole2013spectral}
\[
L=L^{(L)}+L^{(I)}
\]
with the intra-layer-supra Laplacian $L^{(L)}=\mathrm{blkdiag}(L_1,\ldots,L_K)$ and the interlayer-supra Laplacian 
$L^{(I)}=L_I\otimes I$ with $L_I$ the interlayer Laplacian.  Again, the eigeninformation of the supra-Laplacian provides rich information about the network \cite{sole2013spectral,radicchi2014driving,taylor2017eigenvector}. One can also obtain networks where the nodes are not connected across layers but rather only have intra-layer connections (cf. \cite{mercadoNeurIPS2019}). We do not discuss the full details of the various different network structures here but identify this as a very exciting area of future research.

In the context of analyzing social relationships we want to mention \emph{signed networks}, which are graphs with positive and negative edge weights. These networks are used to model friend and foe type relationships and we refer to \cite{leskovec2010signed,sedoc2017semantic,hage1973graph,leskovec2010predicting,tang2016survey} and the references mentioned therein for an overview of some of the crucial applications.
Again techniques such as spectral clustering \cite{sedoc2017semantic,mercado2016clustering,pmlr-v97-mercado19a}, semi-supervised learning \cite{bosch2018node}, convolutional networks \cite{derr2018signed} are available to extract further information from the data. The difficulty for signed networks is that the classical graph Laplacian \change{is not feasible as for example the sum of the weights could be zero resulting in a non-invertible degree matrix. As a result several competing Laplacians are possible (see \cite{gallier2016spectral,tang2016survey} for an overview).}

The graph Laplacian is also an essential tool analyzing complex networks via network motifs \cite{benson2016higher,paranjape2017motifs}, graph centralities \cite{estrada2012physics,estrada2005subgraph,qi2012laplacian}, or community detection \cite{newman2004detecting,newman2006modularity}. It has also been suggested to replace the graph Laplacian by a \emph{deformed Laplacian} or \emph{Bethe Laplacian} \cite{bruna2017community,saade2014spectral,morbidi2013deformed}
as
$$
H(s)=(s^2-1)I-sW+D
$$
for a parameter $s\in\R$ and $W$ the adjacency matrix of the graph. \change{It has been shown that $H(s)$ corresponds to a non-backtracking random walk, which is a simple random walk that is conditioned not to jump back along
the edge it has just traversed \cite{krzakala2013spectral} and shows better performance when community detection is desired in very sparse graphs generated by the stochastic block model.}

In the next section, we turn our attention to the case when the kernel, i.e. \eqref{eq::Gaussian} is not only part of the graph Laplacian, but is viewed as the defining element of embedding the data into a high-dimensional space.
\subsection{Kernel methods}
\begin{figure}[htb!]
	\centering
	\includegraphics[width=.75\linewidth]{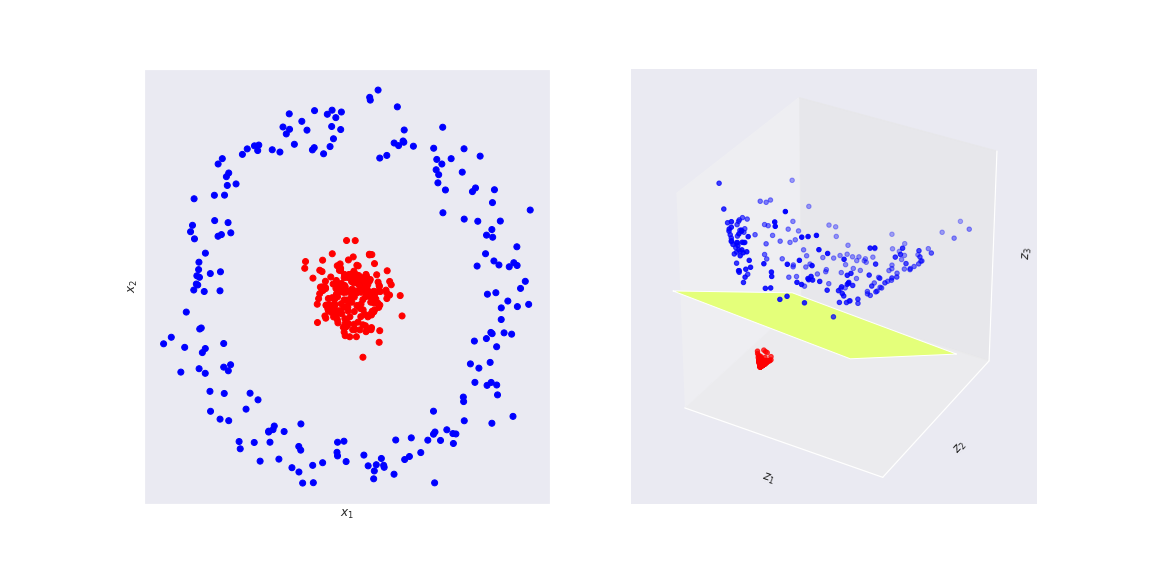}
	\caption{A typical two-dimensional dataset (left) difficult to separate linearly vs the embedding of the data into three-dimensional space (right) with the decision boundary shown in grey.}
	\label{fig:kernel}
\end{figure}

The transformation of the data via a so-called \emph{kernel function}, like the Gaussian encountered for the graph Laplacian, is a technique that has been successfully applied in many data science tasks \cite{shawe2004kernel,scholkopf2001learning,muller2001introduction,hofmann2008kernel}. In \cref{fig:kernel} we illustrate a dataset that is difficult to linearly separate in two dimensions on the left. When the data is transformed via kernelization to three dimensional space it is easily separable.
\change{
Kernel methods can be motivated from the evaluation of the function 
$f(x)=w^Tx$ for a new data point $x$ by using the computed weights $w$ as minimizers of a least squares problem \eqref{ls1}. By employing Lagrangian duality \cite{scholkopf2001learning} we can write the weights as
$w=\frac{1}{\lambda}\sum_i\alpha_i\x_i\in\R^m$  where $\lambda$ is a regularization parameter, $\x_i\in\R^m$ are the vectors associated with the data and the $\alpha_i$ are the Lagrange multipliers.  Inserting the new point $x$ gives 
$
f(x)=\frac{1}{\lambda}\sum_i\alpha_i\x_i^Tx,
$
which relies on the evaluation of inner products $\x_i^Tx$. It turns out that the evaluation of this and other inner products is replaced by a more general \emph{kernel function} $k(\cdot,\cdot)$. We then write the above as 
$$
f(x)=\frac{1}{\lambda}\sum_i\alpha_ik(\x_i,x)
$$
using the kernel instead of the inner product. The goal within kernel methods is now to find a kernel that allows for better separability then the trivial kernel $k(\x_i,x)=\x_i^Tx.$
}
To understand the role of the kernel function let us look at the kernel \cite{scholkopf2001learning}
$$
k(v,x)=(x^Tv)^2$$
and consider the feature mapping $\phi$ that maps the two-dimensional data into three-dimensional space via
$$
\phi(x)=\phi(\begin{bmatrix}
x_1\\x_2\\
\end{bmatrix})=
\begin{bmatrix}
x_1^2\\
\sqrt{2}\ x_1x_2\\
x_2^2
\end{bmatrix}
$$
and we can see that 
$$
k(x,v)=\phi(x)^{T}\phi(v).
$$
This even comes at a computational advantage as we do not need to compute the higher-dimensional $\phi$ vectors. This technique of  avoiding the direct computation of the higher-dimensional nonlinear relation via the evaluation of a kernel function is known as the \emph{kernel trick} \cite{schaback2006kernel,scholkopf2001learning,scholkopf2001kernel}. It is clear now that the adjacency matrix of the graph Laplacian is also a kernel matrix for the particular choice of the Gaussian kernel. 
The assembly of the kernel for all data points $\x_i$ into a matrix leads to a positive semi-definite Gram/kernel matrix $K$. In fact, under the name of \emph{kernel PCA} the leading (now largest) eigenvectors of $K$ are computed to obtain principal components/directions in the data (cf. \cite{scholkopf1998nonlinear,scholkopf1997kernel}). 
The use of the kernel trick in machine learning is omnipresent. One of the simplest but powerful methods is the so-called kernel ridge regression (KRR), where the core problem is to minimize the function
$$
\frac{1}{2}\norm{b-\X u}+\frac{\lambda}{2}\norm{u}
$$
with $b$ a vector encoding the training data and $u$ a vector of weights. \change{This problem is then solved using a dual formulation resulting in a formulation based on the matrix $\X\X^T$, which consists of inner products of the feature vectors. A kernelization of $\X\X^T$ then leads to the kernel matrix $K$ for which we need to solve the linear system \cite{alfke2018nfft}}
\[
\left(K+\lambda I\right)w=b.
\]
The system matrix is symmetric and positive definite for a positive regularization or ridge parameter $\lambda\in\R$. Such a system is typically solved numerically with the use of a preconditioned iterative solver such as the conjugate gradient method \cite{GolVLoan13,saad2003iterative,hestenes1952methods}. The key ingredients are the matrix vector products with $K+\lambda I$ and developing a preconditioner $P\approx K+\lambda I$. For the matrix vector product the non-equispaced fast Fourier transform (NFFT) can be used for a variety of different kernels \cite{alfke2018nfft} and in the case of Gaussian kernels many bespoke methods exist such as \cite{yang2005efficient,yu2017geometry,march2015kernel}. A method based on sketching and the random feature method \cite{rahimi2008random} is introduced in \cite{avron2017faster}, which can also be applied to various different kernel functions. A flexible preconditioner, which changes in every iteration, combined with a suitable Krylov solver was presented in \cite{srinivasan2014preconditioned} whereas the authors in \cite{shabat2019fast} construct a low-rank approximation preconditioner based on an interpolative decomposition and fast matrix vector products. A more difficult scenario arises when $m,$ the dimensionality of the feature vectors $\x_i\R^m$ gets larger. In this case, fast matrix vector multiplication becomes more difficult and suffers from the curse of dimensionality unless certain decay rates are imposed on the matrix entries \cite{morariu2009automatic}. In \cite{you2018accurate} the authors consider an example with $m=90$ where they use a divide-and-conquer parallel algorithm that also comes with communication avoidance. For more details we refer to \cite{you2018accurate} and also the mentioned literature for competing methods. In \cite{rudi2017falkon} the authors use a technique based on randomization (cf. Section \ref{sec::random}) for both the matrix vector product and the preconditioner in KRR while a similar problem is solved in \cite{chenhan2016inv} via high performance computing approaches.

The power of the kernelization has been exploited as an essential ingredient of support vector machines (SVM) \cite{scholkopf2001learning,Vapnik:1982:EDB:1098680,cortes1995support}. \change{In more detail, support vector machines are derived from maximizing the margin of the separating hyperplane. The support vectors are the closest data points to this hyperplane. In order to be able to obtain a nonlinear hyperplane kernelization of the inner products is employed. }The resulting kernel matrix is then at the heart of the quadratic program that needs to be solved for determining the support vectors \cite{scholkopf2001learning} but the computation suffers from having to deal with dense kernel matrices. In \cite{platt1998sequential} the sequential minimal optimization (SMO) method is introduced, which breaks the problem into smaller, analytically solvable problems. Kernel matrices are also crucial when solving the closely related  support vector regression (SVR) problem \cite{smola2004tutorial,drucker1997support}. While both SVM and SVR remain popular methods, deep learning techniques have recently gained more popularity and have been combined with  kernel techniques via graph convolution networks \cite{henaff2015deep,kipf2016semi} or as loss functions in convolutional networks \cite{tang2013deep,madjarov2012extensive}.

Often the linear algebra tasks associated with large scale kernel methods will involve randomization, which we discuss next.
\section{Randomization}
\label{sec::random}
\change{
In their seminal review paper \cite{halko2011finding} the authors discuss the importance of the decompositional approach to matrix computations \cite{dongarra2000guest} (cf. Section \ref{sec::fac}).  Due to the computational complexity it is not always straightforward to efficiently compute matrix factorizations, such as the SVD or QR, since in many data science applications the matrix $\X$ is so large that computing approximate factorizations is too costly.} Additionally, the data might be corrupted or it might be desirable to avoid many passes over the data so that classical methods need further thought. 

One of the key ingredients in allowing efficient numerical linear algebra is the process of \emph{randomization.} Randomization has proven to be a valuable tool in various matrix computation tasks such as matrix vector products \cite{cohen1999approximating,drineas2006fast} or low-rank approximations \cite{bingham2001random,drineas2006fast2}. These techniques typically follow the scheme of producing a skinny matrix $Q\in\R^{n,k+p},$ \change{with $k$ the desired approximation rank and $p$ an oversampling parameter that is needed for theoretical guarantees,} such that we obtain the orthogonal-projection-approximation onto the subspace spanned by $Q$ via
$$
A\approx QQ^{T}A.
$$
From the information contained in $Q$ one then computes standard factorizations such as QR or SVD decompositions. 
In more detail, we form the matrix $B=Q^{T}A\in\R^{k+p,m}$, which then gives the low-rank approximation $A\approx QB.$ \change{ Replacing $B$ by its SVD results in an approximate SVD of $A$ via
$$
A\approx QB=Q(U\Sigma V^T)=(QU)\Sigma  V^T.
$$ }
The first stage, i.e, the computation of the matrix $Q$, follows from the prototype algorithm illustrated in Algorithm \ref{alg::pro1}. 
\begin{algorithm}
	\caption{Prototype algorithm from \cite{halko2011finding}.}
	\label{alg::pro1}
	\begin{algorithmic}
		\State Draw a random $n \times (k + p)$ test matrix $G$\footnotemark.
		\State Compute $Y=AG$.
		\State Construct orthogonal basis for $\mathrm{range}(Y),$ e.g. $[Q,R]=qr(Y)$.
	\end{algorithmic}
\end{algorithm}
\footnotetext{Typically, this is drawn as a Gaussian random matrix.}
We are typically only interested in a decomposition of order $k$ when proceeding to the second stage of the method but that $G$ is of dimension $k+p$ where $p$ is the oversampling parameter. The dimensionality of $G$ is crucial as our aim is to reduce the complexity as much as possible. Theoretical bounds for the approximation quality are given in \cite{halko2011finding}  and depend on the parameters $k$ and $p$ as well as the matrix size and the $(k+1)$-st singular value of $A$. Such a randomized SVD has become the an essential tool in many scientific computing and data science applications in areas such as uncertainty quantification \cite{li2019sparse}, optimal experimental design \cite{alexanderian2014optimal}, or computer vision \cite{erichson2016randomized}.

One of the key applications of randomized methods within data science is the use for approximation of kernel matrices \cite{martinsson2018randomized,drineas2005nystrom,zhang2008improved,li2014large,hsieh2014fast, bertozzi2012diffuse} where the desire is to avoid the computation of the full kernel matrix since the storage demand can be too large. In particular, the kernel matrix
is approximated via 
$$
A\approx 
(AQ)(Q^TAQ)^{-1}Q^TA^T
$$
where in the traditional setup the matrix $Q$ contains columns of the identity matrix drawing columns from the matrix $A$
and we obtain the so-called \emph{Nystr\"om scheme}
$$
\begin{bmatrix}
A_{11}\\
A_{21}\\
\end{bmatrix}
A_{11}^{-1}
\begin{bmatrix}
A_{11}&A_{12}\\
\end{bmatrix}
=
\begin{bmatrix}
I\\
A_{21}A_{11}^{-1}\\
\end{bmatrix}
\begin{bmatrix}
A_{11}&A_{12}\\
\end{bmatrix}
=
\begin{bmatrix}
A_{11}&A_{12}\\
A_{21}&A_{21}A_{11}^{-1}A_{12}\\
\end{bmatrix}
\approx A.
$$
This means in order to approximate a kernel matrix $A$ one only needs to compute a small number of columns and rows, which has been shown to be very successful (cf. \cite{pourkamali2018randomized,drineas2005nystrom}). \change{The authors in \cite{drineas2005nystrom} give an approximation bound using the diagonal entries of the kernel matrix and a probabilistic interpretation whereas the authors in \cite{osinsky2018pseudo} provide a bound for the Chebyshev-norm of the pseudo-skeleton approximation and the $k+1$-st singular value of the kernel matrix.} In \cite{martinsson2018randomized} the author suggests to use a $Q$ obtained via Algorithm \ref{alg::pro1} and thus obtain a variant of the traditional Nystr\"om method.

We already pointed out in Section \ref{sec::ID} that interpretable decompositions can be obtained via QR or SVD approximations to $\X$. As the author in \cite{martinsson2018randomized} points out computing a matrix that holds a basis for the column space of $A$ is crucial. For this the author in \cite{martinsson2018randomized} computes $Y=\X G$ and via a subspace iteration proceeds to iteratively refine the matrix $Y$ using the update $Y\leftarrow AA^{T}Y$. The resulting $Y$ is then used to identify the relevant row-indices for the ID. In the same way one can obtain the two-sided ID and also a randomized CUR decomposition, where the details are given in \cite{martinsson2018randomized}. Note that algorithms for directly computing the CUR decomposition compute some sort of statistical score for the importance of the columns/rows via the truncated SVD \cite{mahoney2009cur,wang2013improving,boutsidis2009improved}. When obtaining the CUR from the two-sided ID we get the index sets from the ID matrices. The usefulness of the importance sample or sampling statistics is illustrated in \cite{mahoney2011randomized}  where the author shows that operations from numerical linear algebra such as inner products or matrix vector products can be performed using randomized numerical linear algebra. The key ingredients when multiplying $AB$ or $a^Tb$ is the sampling strategy for producing the approximations. The strategy is often based on an \emph{importance sampling} strategy \cite{drineas2012fast,drineas2016randnla} to draw elements or rows/columns of a matrix for further processing. For example, the authors in \cite{eriksson2011importance} use the vector $q$ 
$$
q_j=\frac{\norm{A(:,j)}_2\norm{B(j,:)}_2}{\sum_{i=1}^{n}\norm{A(:,i)}_2\norm{B(i,:)}_2},\quad j=1,\ldots,n
$$
of probabilities\footnote{This is obtained from the minimization of the expected value of the variance of the inner product/matrix vector product.} 
to produce an approximation $C\approx AB$ based on rows of $A$ and columns of $B$. For the solution of a least squares problem such as \eqref{ls1} \cite{drineas2011faster,mahoney2011randomized} the (random) selection process is decoupled from from a deterministic computation via a traditional numerical linear algebra method. In particular, the first stage computes a score based on e.g. Euclidean norms, a truncated SVD matrix, or columns of a truncated Hadamard matrix. This means that the original problem is separated into a random sampling procedure leading to a reduced formulation that is then solved via a deterministic NLA approach.
Due to the success in data science applications randomized methods have also penetrated classical problems in scientific computing such as solving linear systems of equations \cite{gower2015randomized,tu2017breaking,needell2010randomized}, eigenvalue problems \cite{saibaba2016randomized,gu2015subspace} or inverse problems \cite{xiang2015randomized,xiang2013regularization,vatankhah2018total}. A recent survey can be found in \cite{martinsson2020randomized}.

One area where randomization has helped greatly with the evaluation of complex mathematical expressions is the computations of functions of matrices, which we want to describe now.

\section{Functions of matrices}
\label{sec::func}
The concept of evaluating a function of a matrix 
$$
f:R^{n,n}\rightarrow R^{n,n},
$$
where $f$ is not meant to be evaluated element-wise, is an old but still very relevant one and we refer to \cite{higham2008functions} for an excellent introduction to the topic.
\change{
Matrix functions appear in a variety of applications with a particularly important example given by Gaussian processes \cite{rue2005gaussian}; an ubiquitous tool in machine learning and statistics or in the analysis of complex networks \cite{estrada2010network,benzi2020matrix}.}
One can define the matrix function via the Cauchy integral theorem as
\begin{equation}
\label{eq::cauchy}
f(A)=\frac{1}{2\pi i}\int_\Gamma f(z)(zI-A)^{-1}dz
\end{equation}
provided $f$ is analytic on and inside a closed contour $\Gamma\subset\mathbb{C}$ that encloses the spectrum of $A$. There are other definitions that are equivalent for analytic functions (cf. \cite{higham2008functions}). Matrix functions provide a large number of challenges for numerical methods and an early discussion with 19 \textit{dubious methods} for the matrix exponential can be found in \cite{moler1978nineteen} with an update provided in \cite{moler2003nineteen} adding Krylov subspace methods as the twentieth method. One typically considers two different setups, the first is the computation or approximation of the matrix $f(A)$ explicitly and the other is the evaluation of $f(A)b$ where the explicit computation of $f(A)$ and subsequent application to $b$ is avoided. For a recent survey on evaluating $f(A)b$ we refer to \cite{GuettelKressnerLund2020}. In this case the computation of $f(A)$ and then applying it to $b$ is obviously avoided. Efficient techniques typically depend on the size of the data matrix $A$. The approximation of $f(A)b$ employing the Cauchy integral formula \eqref{eq::cauchy} via contour integration was given in \cite{hale2008computing}. Often the main work goes into solving linear systems with $(zI-A)$ and once direct solvers become infeasible approaches based on Krylov subspace methods have shown very good performance for approximating the evaluation of $f(A)b$ \cite{hochbruck1997krylov,knizhnerman2010new,druskin1998extended,guttel2013rational,knizhnerman2010new}. For the more complex case of computing $f(A)$ we again refer to \cite{higham2008functions} for an overview of suitable methods. We here want to mention certain scenarios in data science applications that lead us to such matrix functions.

Throughout scientific computing and engineering fractional Laplacians, where instead of second order derivatives one considers derivatives of arbitrary orders, have gained much popularity in recent years \cite{podlubny1998fractional} and one way to evaluate this efficiently is using matrix functions, e.g. \cite{burrage2012efficient}. Fractional powers of the Laplacian have recently become more popular for semi-supervised learning \cite{bautista2019gamma,de2017fractional}, where to the best of our knowledge techniques based on the full eigendecomposition have been employed, which for large graphs quickly becomes infeasible. In \cite{mercado2018power} a multilayer graph is considered where arbitrary matrix powers of the layer Laplacians are evaluated using the polynomial Krylov subspace method (PKSM) \cite{higham2008functions} and using contour integrals \cite{mercadoNeurIPS2019}. 

We return again to the study of complex networks and are now interested in the network centrality \cite{estrada2010network} \change{as this helps in identifying the most important nodes, such as the most influential people in a social network. We assume that the network is understood as an undirected graph and describe some of the most important centrality measures.} One measure of centrality is the degree matrix $D$ of the graph Laplacian assigning the degree $d_{ii}$ to every node. Other centrality measures are of great use for understanding the underlying data. For example, the eigenvector centrality \cite{bonacich1987power} is defined as
$$
b_i=\frac{1}{\lambda_1}W\phi_1
$$
with $W$ the adjacency matrix of the network and $(\lambda_1,\phi_1)$ the Perron-Frobenius eigenvalue and eigenvector, respectively. The authors in \cite{estrada2010network} illustrate that powers of the adjacency matrix $W$ provide essential information about the graph. In particular, studying the $(i,j)$-th entry of $W^n$ allows to count the number of different walks of length $n$ getting from node $i$ to node $j$. For example the diagonal entries of $W^2$ give the degrees of the nodes in the network. Note that powers of the graph Laplacian are considered in \cite{abbe2018graph} for the purpose of spectral clustering and within graph convolutional networks for semi-supervised learning in \cite{liu2019higher}. The use of $W^2$ suggests to consider higher powers of the adjacency matrix for longer walks and one obtains
$$
(I+\frac{W^2}{2!}+\frac{W^3}{3!}+\frac{W^4}{4!}+\ldots)_{ii}=\exp(W)_{ii}=e_i^T\exp(W)e_i
$$
as the so-called subgraph centrality for node $i$, (cf. \cite{estrada2000characterization,estrada2005subgraph}) with $e_i$ the corresponding unit vector. The \change{quantity $\mathrm{trace}(\exp(W))$} is defined as the Estrada index of a network \cite{ginosar2008estrada,deng2009estrada,de2007estimating}. Due to the high complexity the computation of the matrix function $\exp(W)$ should be avoided, especially if complex networks are considered. \change{It is well known that expressions of the form
$
u^Tf(A)u
$ can be well approximated by $e_1^Tf(T_k)e_1,$ where $T_k$ is the tridiagonal matrix coming from the Lanczos process applied to $A$ and $u$ is assumed to have norm one. In particular, for our example we have $u=e_i$, $f(\cdot)=\exp(\cdot)$, and $A=W$.} As the authors in \cite{golub2009matrices,golub1994matrices} show there is a beautiful relation between $
u^Tf(A)u
$ and Gauss--quadrature with first results going back to \cite{golub1969calculation}. In particular, the authors in \cite{benzi2010quadrature} use the error estimates that stem from studying the relation between the Lanczos process and Gauss quadrature. In \cite{bonchi2012fast} the authors also rely on the power of Gauss quadrature to compute Katz scores and commute times between nodes in a network.
\change{
Using the Lanczos process or other Krylov subspace methods \cite{golub2009matrices,golub2008approximation,strakovs2011efficient} to approximate quantities of the form $u^Tf(A)u$ also proves essential for many applications (cf. \cite{avron2010counting,hutchinson1990stochastic,atallah2001randomized}) when the trace of a matrix (function) has to be estimated via
$$
\mathrm{trace}(f(A))\approx \frac{1}{s}\sum_{i=1}^{s}v_i^{T}f(A)v_i
$$
where we rely on $s$ carefully chosen random vectors $v_i$. One such estimator is the so-called Hutchinson estimator \cite{avron2011randomized} where $v_i=\pm 1$ and the numerical computation again relies on the 
Lanczos process and its efficient approximation of the spectrum of $A$ \cite{ubaru2017fast,meurant2009estimates, bellalij2015bounding}.} 
\change{
Gaussian processes \cite{rasmussen2003gaussian,mackay1998introduction,rue2005gaussian} are another essential tool within statistics and data science. Their evaluation poses a significant numerical challenge as it relies on evaluating matrix functions as we will illustrate now.} Following \cite{dong2017scalable}, a Gaussian process is a collection of random variables with a joint probability distribution. Let us consider $X=\left\lbrace x_1,\ldots,x_n\right\rbrace$ with all the $x_i\in\R^d$. 
The Gaussian process can then define a distribution over functions 
$f(x)\sim\mathcal{GP}(\mu(x),k(x,x'))$ with mean $\mu(x):\R^d\rightarrow \R$ and covariance function $k(x,x'):\R^d\times \R^d\rightarrow\R$. Now $f_X\in\R^n$ represents the vector of the function values for $f(x_i)$, $\mu_X\in\R^n$ the evaluation of $\mu(x_i)$, and the matrix $K_{XX}\in\R^{n,n}$ represents the evaluation of $k(x_i,x_j)\quad\forall i,j$. Then it holds that $f_X\sim\mathcal{N}(\mu_X,K_{XX}).$ The chosen covariance kernel depends on hyperparameters $\theta$ and we assume that we are given a set of noisy function values $y\in\R^n$ with variance $\sigma^2.$ Assuming a Gaussian prior distribution depending on $\theta$ one obtains the log marginal likelihood as
$$
\mathcal{L}(\theta|y)=-\frac{1}{2}\left[\left(y-\mu_X\right)^T\left(K_{XX}+\sigma^2I\right)^{-1}\left(y-\mu_X\right)+\mathrm{log}\left(\mathrm{det}\left(K_{XX}+\sigma^2I\right)\right)+n\mathrm{log}(2\pi)\right].
$$
It is obvious that for the numerical solution of the minimization of the log marginal likelihood the evaluation of the log-determinant is crucial. If the matrix $K_{XX}$ is small then we can afford the computation via the Cholesky decomposition, which is an $\mathcal{O}(n^3)$ computation. An efficient computation is then given via
$$
\mathrm{log}\left(\mathrm{det}(W)\right)=2\sum_{i=1}^{n}L_{ii}
$$
where $W=LL^T$ is the Cholesky decomposition of the symmetric positive definite matrix $W$. In  \cite{rue2001fast, rue2005gaussian} the authors show that the methods for approximating the log-determinant become more efficient if one works with Markovian fields (cf. \cite{simpson2008fast} for approaches beyond Cholesky for Markovian fields). In general the Cholesky decomposition will be too expensive and again Krylov subspace methods such as the Lanczos process and its connection to Gauss--quadrature are exploited. From the relation
$$
\exp(\mathrm{trace}(\mathrm{log}(W)))=\exp(\mathrm{trace}(X\mathrm{log}(\Lambda)X^{-1})=\exp\left(\sum_i\mathrm{log}(\lambda_i)\right)=\prod_{i}^{}\lambda_i=\mathrm{det}(W)
$$
we obtain the following beautiful relationship
$$
\mathrm{log}(\det (W)) = \mathrm{trace}(\mathrm{log}(W)).
$$
\change{
This shows that in order to approximate the log--determinant, it is again crucial to efficiently evaluate a matrix function, namely 
$$
\mathrm{trace}(\mathrm{log}(A))\approx \frac{1}{s}\sum_{i=1}^{s}v_i^{T}\mathrm{log}(A)v_i.
$$} The Lanczos-based approximation has been used in \cite{lin2016approximating,bai1998computing,ubaru2017fast} and the authors in \cite{dong2017scalable} combine these ideas with fast matrix vector products by modified structured kernel interpolation \cite{wilson2015kernel} allowing a good approximation of the log--determinant and also its derivatives.

We refer to \cite{higham2016catalogue} for an overview of software for the computation of matrix functions in various programming languages.

\section{High-dimensional problems}
\label{sec::tensors}
All of the above problems are tailored to when the data is given as a matrix $\X\in\R^{n,m}$ but often it is also natural for the data to be given as a tensor 
\begin{equation}
\X\in \R^{n_1\times n_2\times \ldots\times n_d}
\label{eq::tensor}
\end{equation}
with the dimensions $n_i\in\N$ for all $d$ modes. It is obvious that the storage requirements for a full tensor $\X$ quickly surpass the available resources in many applications. This exponential increase in relation to the parameter $d$ is often referred to as the \emph{curse of dimensionality.} Tensors and their approximations are closely related to the representation of a nonlinear function in a high--dimensional space \cite{devore1998nonlinear}
\begin{equation}
f(x_1,\ldots,x_d)\approx \sum_{j=1}^{l}\alpha_j\phi(x_1,\ldots,x_d)
\label{eq::hdfunc}
\end{equation}
where the choice of functions $\phi$ is crucial and these could be a basis, a frame, or a dictionary of functions \cite{temlyakov2011greedy,bungartz2004sparse,cohen2015approximation,dung2018hyperbolic,kammerer2015approximation}. In this case it is desirable to find an approximation to $f(x_1,\ldots,x_d)$ with some form of sparsity, i.e., as small as possible number of terms $l$. Alternatively, one can interpret the function as being defined on a tensor product space and aiming at designing a low-rank approximation to $f$ (cf. \cite{nouy2019higher} for example). The approximation of $f$ in \eqref{eq::hdfunc} is then performed by relying on similar formats and tool as the ones that are used for the approximation of the given tensor $\X$ in \eqref{eq::tensor}, which we briefly describe now.

For this reason approximating the tensor has a long tradition in computational mathematics, physics, chemistry, and other disciplines. While being already 10 years old we would like to point to \cite{kolda2009tensor} for a seminal overview of numerical tensor methods \change{due to the introductory character and the broad overview.} There are several low-rank approximations to the full tensor $\X$ that correspond to different tensor formats such as the classical CANDECOMP/PARAFAC(CP) formulation \cite{kiers2000towards} where
\begin{equation}
\X=\sum_{j=1}^{R}\lambda_j u_{j}^{(1)}\circ\cdots \circ u_{j}^{(d)},
\end{equation}
with $\circ$ being the dyadic product. While this format is the most natural, it has properties that make it not ideal for (many) applications, e.g. the set of rank $R$ tensors is not closed, the computation ill-posed, etc. (cf. \cite{grasedyck2013literature} for more details). Alternatively, one can consider the Tucker format/HOSVD \cite{de1997signal, de2000multilinear}, which somewhat generalizes the SVD to the multidimensional case but still suffers from the curse of dimensionality due to the $d$-dimensional core tensor needed. This is overcome by the tensor-train format \cite{oseledets2011tensor}, which consists of $d$ core tensors of maximal dimension $3.$ In the context of data science tensors have long been an essential tool and to list all applications and techniques goes beyond the scope of this paper.  We here only list a few results that address similar questions to the above mentioned matrix cases: non-negative factorizations \cite{cichocki2007nonnegative,cichocki2009nonnegative,chi2012tensors}, interpolative decomposition \cite{saibaba2016hoid}, CUR \cite{cortinovis2019low,mahoney2008tensor}, kernel methods \cite{tao2005supervised,he2017kernelized,he2014dusk}, semi-supervised learning \cite{gujral2018smacd}, and randomization \cite{battaglino2018practical,wang2015fast}.

This list is far from complete and for more detailed reviews on tensor methods for data science we refer to the recent surveys \cite{cichocki2016tensor,cichocki2017tensor,cichocki2014tensor} and the references given therein. For a very broad overview of recent techniques and literature for numerical tensor methods beyond data science we refer to  \cite{grasedyck2013literature}. 

One of the key areas in which tensor prove very powerful tools is the field of deep learning, which we now briefly review.
\section{Numerical linear algebra in deep learning}
The study of computational and mathematical aspects of deep learning is glowing white hot with activity and we only want to point to number of places where linear algebra is used for better understanding or improving performance of deep learning methods. 

In the most generic setup, the task in an artificial neural network is to determine the weight matrices $W^{(j)}$ and bias vectors $b^{(j)}$ from minimizing a loss function to obtain the function
\begin{equation}
F\left(x\right):=W^{\left(L\right)} \left(\ldots\sigma \left(W^{\left(2\right)} \left(\sigma\left(W^{\left(1\right)}x+b^{\left(1\right)}\right)\right)+b^{\left(2\right)}\right)\right)\ldots +b^{\left(L\right)}
\end{equation}	
see \cite{higham2019deep, goodfellow2016deep,lecun2015deep} for more information and further references. The loss function typically consists of a sum of many terms due to the large number of training data and as a result the optimization is based on gradient descent schemes \cite{bottou2010large}. For classical optimization problems it is often desirable to use second order methods due their superior convergence properties. These typically require the use of direct or iterative solvers to compute the step towards optimality, e.g. solving with the Hessian matrix in a Newton method. As a result the applicability of Newton-type schemes for the computation of $W^{(j)}$ and $b^{(j)}$ has received more attention with a strong focus on exploiting the structure of the Hessian matrix \cite{botev2017practical,dangel2019modular,wang2018newton,chen2019fast,sagun2016singularity,vinyals2012krylov}. 

The computational complexity of neural networks is challenging on many levels. In order to reduce the numerical effort the authors in \cite{teoh2006estimating,xue2013restructuring} use the SVD to analyze the weight matrices and then restructure the neural network accordingly, while the authors in \cite{yousefzadeh2019refining} use the condition numbers of the weight matrices for such a restructuring. The weight matrices have received further attention as in \cite{novikov2015tensorizing} the authors compress the weight matrices of fully connected layers using a low-rank tensor approximation, namely, the tensor train format \cite{oseledets2011tensor}. Additionally, low-rank approximations are also useful for speeding up the evaluation of convolutional neural networks \cite{jaderberg2014speeding} by using a low--rank representation of the filters, which are used for detecting image features. For a very similar task the authors in \cite{lebedev2014speeding,gusak2019automated,gusak2019musco} rely on optimized tensor decompositions. Note that recently these techniques have also been applied to adversarial networks \cite{cao2018tensorizing}.

The connection of neural network architectures to optimal control problems as introduced in \cite{haber2017stable,gunther2018layer} makes heavy use of PDE-constrained optimization techniques including efficient matrix vector products. This topic has also recently received more attention from within the machine learning community \cite{chen2018neural} and promises to be a very interesting field for combining traditional methods from numerical analysis with deep learning.

As already pointed out in Section \ref{sec::GL} graph convolutional networks (GCNs) have shown great potential and received much popularity recently as tools for semi-supervised learning \cite{kipf2016semi}. The expressive power of representing data as graphs that is the basis of \change{GCNs has led to using their methodology for adversarial networks \cite{DBLP:journals/corr/abs-1906-03220}, matrix completion \cite{berg2017graph}, or within graph auto encoders \cite{kipf2016variational}.}

This is only a brief glimpse where numerical linear algebra techniques have enhanced deep learning methodology. We believe that this will be an area of intense research in the future.
\section{Conclusion}
We have illustrated with this literature review that numerical linear algebra is alive and kicking. In a sense the rise of machine learning, big data, and data science shows that \textit{old methods} are still very much in fashion but also that new techniques are required to either meet the demands of practitioners or account for the complexity and size of the data.
\section*{Acknowledgments} The author would like to thank Dominik Alfke, Peter Benner, Pedro Mercado, and Daniel Potts for helpful comments on an earlier version of this manuscript. He is also greatly indebted to the two anonymous referees who greatly helped to improve the presentation.

\bibliographystyle{siam}
\bibliography{nlads_refs}
\end{document}